\documentclass{article}
\usepackage{amsfonts}
\usepackage{amsfonts}
\usepackage{amssymb}
\newcommand{\lbl}[1]{\label{#1}}

\newcommand{\rf}[1]{(\ref{#1})}
\normalbaselines
       \newtheorem{Theorem}{Theorem}[section]
       \newtheorem{Proposition}{Proposition}[section]
       
       \newtheorem{Example}{Example}[section]
       \newtheorem{Lemma}{Lemma}[section]
       \newtheorem{Remark}{Remark}[section]
 \newcommand{\QED}{{\fbox{}\newline }}

\newcommand{\ZZ}{ Z\!\!\!Z }
\newcommand{\NN}{I\!\!N}

\newcommand{\be}{\begin{equation}}
\newcommand{\ee}{\end{equation}}
\newcommand{\bea}{\begin{eqnarray}}
\newcommand{\eea}{\end{eqnarray}}
\newcommand{\conj}{\widetilde{=}}

\newcommand{\qed}{\QED}

\newenvironment{proof}{\noindent {\bf Proof.\/}}{\qed\vskip 0.1in}
\newenvironment{proofof}[1]{\noindent {\bf Proof of #1.\/}}{\qed\vskip 0.1in}
\author{Wlodzimierz  Bryc \\ 
 Department of Mathematics \\
University of Cincinnati\\
PO Box 210025\\
Cincinnati, OH 45221--0025\\
Wlodzimierz.Bryc@UC.edu \thanks{ 
\newline
{\bf AMS (1991) Subject Classification:}  60E99}
}
\title{Stationary Markov chains with linear regressions}
\begin{document}
\maketitle
\begin{abstract} In  
Bryc(1998)
 we determined one
dimensional distributions of a stationary field with linear regressions
(\ref{(A')}) and quadratic conditional variances  (\ref{(B')}) under a
linear  constraint (\ref{MainAssumption}) on  the coefficients of the
quadratic expression (\ref{Q}). In this paper we show that  for stationary
Markov
chains with linear regressions and quadratic conditional variances 
  the coefficients of the
quadratic expression are indeed tied by a linear constraint which
 can take only one of the two alternative forms  (\ref{MainAssumption}),  or
(\ref{degenerate1}).

\end{abstract}

\section{Introduction}

Let $(X_k)_{k\in \ZZ}$ be  a  square-integrable random sequence.
Consider the following two conditions.
\be \lbl{(A')}
E(X_k| \dots,
X_{k-2},X_{k-1},X_{k+1},X_{k+2},\dots)=L(X_{k-1}, X_{k+1})  
\ee
for all $k\in\ZZ$.
\be \lbl{(B')}
E(X_k^2| \dots,
X_{k-2},X_{k-1},X_{k+1},X_{k+2},\dots)=Q(X_{k-1}, X_{k+1})
\ee
for all $k\in\ZZ$.

A number of  papers 
analyzed
conditions similar to (\ref{(A')}) and 
(\ref{(B')}). Of particular interest are papers   
Wesolowski(1989) and Wesolowski(1993),  
who analyzed continuous time
 processes $X_t$ with linear regressions and quadratic second order conditional moments $Q()$
 under the assumption that variances of $X_t$ are strictly increasing;
 these processes turned out to have independent increments.
Szablowski(1989) relates distributions of
mean-square differentiable processes to conditional variances.
 Bryc~\&~Plucinska(1983) show that 
 linear regressions and constant conditional variances characterize
 gaussian sequences.  In Bryc(1998)
we show that a certain class of  quadratic functions
$Q$
determines the univariate distributions  for stationary processes which satisfy
(\ref{(A')}) and  (\ref{(B')}) with linear $L$. 
For additional  references the reader is referred to Bryc(1995).

In this paper we assume that $(X_k)$ is strictly stationary and the regressions are given by
 a symmetric linear polynomial $L(x,y)=a(x+y)+b$, 
and a general symmetric quadratic
 polynomial 
\be \lbl{Q}
Q(x,y)=A(x^2+y^2)+B xy + C +D(x+y)\ee
The linear polynomial $L()$ is determined uniquely by the
 covariances of $(X_k)$. Namely, if the random variables $X_k$ 
are centered with variance 
$1$, the correlation coefficients $r_k=corr(X_0,X_k)$, and $r_2>-1$, then
$L(x,y)=\frac{r_1}{1+r_2}(x+y)$. Since the moments of both sides of (\ref{(B')})
 must match,  after
standardization we also get  the trivial
relation
\be \lbl{C}
C=1-2A - B r_2
\ee
This still leaves three parameters $A,B$, and $D$ undetermined. 

In this paper we analyze in more detail 
 which quadratic polynomials $Q()$ can occur in
(\ref{(B')}) when   $(X_k)$ is a stationary Markov chain.
We show that in this case we necessarily have
$D=0$ and that the remaining two coefficients satisfy one
of the two linear equations (\ref{MainAssumption}) or (\ref{degenerate1}).
We  show that if condition (\ref{MainAssumption}) is satisfied then the
 remaining free coefficient satisfies certain inequalities; under additional assumption
(\ref{(A')}), (\ref{(B')}), and (\ref{MainAssumption})
 characterize certain Markov chains uniquely.

\section{Results}
Through the rest of the paper we assume that $(X_k)$ is standardized, $E(X_k)=0, E(X_k^2)=1$.
We denote the correlations by
 $r_k:=E(X_0X_k)$, $r:=r_1$.

For Markov chains the regression equations (\ref{(A')})
 and (\ref{(B')})
 become respectively
\be \lbl{(A)}
E(X_k| X_{k-1},X_{k+1})=L(X_{k-1}, X_{k+1})  
\ee
\be \lbl{(B)}
E(X_k^2|X_{k-1},X_{k+1})=Q(X_{k-1}, X_{k+1})
\ee

The following result shows that the coefficients of (\ref{Q})
are tied by a linear constraint.\begin{Theorem}\lbl{Corollary}
Let $(X_k)$ be a square-integrable 
standardized stationary homogeneous Markov chain such that $r\ne0$, and $2|r|<1+r_2$. 
If $(X_k)$ satisfies conditions
(\ref{(A)}) and (\ref{(B)}), then the coefficients of $Q()$ in (\ref{Q}) 
satisfy $D=0$ and either 
\be\lbl{MainAssumption}
A(r^2+1/r^2)+B=1
\ee  or 
\be \lbl{degenerate1}
2 A + Br^2=1
\ee
(When  $Q$ is non-unique this should be interpreted  that
there is a quadratic function $Q$
with the coefficients satisfying $D=0$ and at least one of the identities 
(\ref{MainAssumption}) or (\ref{degenerate1}).)
\end{Theorem}

It turns out that 
(\ref{MainAssumption}) implies
 additional restrictions on the range of the remaining
 free parameter $A$. 

\begin{Theorem}\lbl{T4}
Let $(X_k)$ be a standardized strictly stationary square-integrable sequence such that
conditions (\ref{(A')}) and (\ref{(B')}) hold true, and the correlation coefficients satisfy
$r\ne0$, and $2|r|<1+r_2$. Suppose that the coefficients of quadratic form $Q()$ in (\ref{Q}) 
are such that $D=0$ and 
(\ref{MainAssumption}) holds true.

Then either $A\ge 1/(1+r^2)$ or  $ A \leq \frac{r^2}{1+r^4}$.

\end{Theorem}

 The next theorem is a version of Bryc(1998), Theorem 2.1.
\begin{Theorem}\lbl{T+}
Suppose that $(X_k)$ satisfies the assumptions of Theorem \ref{T4}, 
and $\frac{r^2}{(1+r^2)^2}\leq A \leq \frac{r^2}{1+r^4}$. 
Then
$X_k$ is a Markov chain with
uniquely determined distribution.
\end{Theorem}
One can also show that  condition \rf{degenerate1}
implies that $|X_k|=|X_0|$ with probability one.

\section{Two-valued Markov chains}\lbl{TwoValued}
Verification of condition \rf{(A)} for two-valued Markov chains is a simple exercise. 
We include it here because two-valued chains
play  a role in the proofs of Theorem \ref{Corollary} and Proposition \ref{T1}.
They also
 occur as "degenerate cases" in linear
regression problems: in Bryc(1998) 
 we construct Markov chains that satisfy \rf{(A)} and \rf{(B)} for 
 $A<r^2/(1+r^4)$; the boundary value $A=r^2/(1+r^4)$  corresponds to
the two-valued case.

We consider only standardized chains with mean $0$ and variance $1$. 
Under this assumption, if a transition matrix is defined by 
\be \lbl{two1}
\Pr(a,a)=1-\alpha, \Pr(a,b)=\alpha, \Pr(b,a)=\beta
, \Pr(b,b)=1-\beta
\ee
 then the invariant distribution assigns probabilities

\be \lbl{two2}
\mu(a)=\frac{\beta}{\alpha+\beta}, \ \mu(b)=\frac{\alpha}{\alpha+\beta}
\ee
 and the two values of the chain are

\be \lbl{two3}
a=\sqrt{\frac{\alpha}{\beta}}, b=-\sqrt{\frac{\beta}{\alpha}}
\ee

We 
 consider non-degenerate Markov chains with the correlation coefficient
$r\ne 0,\pm 1$ only. This excludes three uninteresting
cases: i.i.d sequences, constant
 sequences with $X_k=X_0$ for all $k$, and alternating sequences
 with $X_k=(-1)^kX_0$ for all $k$. 

\begin{Proposition}\lbl{L1}
If $(X_k)$ is a two-valued stationary Markov chain with the one-step 
correlation coefficient $r\ne 0,\pm 1$ then $(X_k)$ satisfies condition (\ref{(A)}) 
if and only if
 $X_0$ is symmetric with values $\pm 1$.
\end{Proposition}
\begin{proof}
First notice that $\alpha\beta>0$, so the values and probabilities in
(\ref{two2}) and (\ref{two3}) are well defined. Indeed, if $\alpha\beta=0$ then we have
$X_k=X_{k-1}$ and
hence
$r=1$.

A simple computation using (\ref{two1}-\ref{two3}) shows that the one-step
correlation coefficient is $r=1-\alpha-\beta$, 
and the two step correlation is $r_2=r^2$.  Since by assumption $0<|r|<1$,  
this implies that 
$\alpha+\beta<2$ and
$\alpha+\beta \ne 1$.

By routine computation we get the following conditional probabilities

%
\bea
\Pr(X_k=a|X_{k-1}=a,X_{k+1}=b)&=&\frac{1-\alpha}{2-\alpha-\beta}
\nonumber\\
\Pr(X_k=b|X_{k-1}=a,X_{k+1}=b)&=&\frac{1-\beta}{2-\alpha-\beta}
\nonumber
\eea

Using (\ref{(A)}) we have
$E(X_1|X_0=a,X_2=b)=\frac{r}{1+r^2}(X_0+X_2)=
\frac{\alpha-\beta}{\sqrt{\alpha\beta}}\frac{1-\alpha-\beta}{1+(1-\alpha -\beta)^2}$. On the
other hand, direct computation using conditional
probabilities
gives
$E(X_1|X_0=a,X_2=b)=
\frac{\alpha-\beta}{\sqrt{\alpha\beta}}\frac{1-\alpha-\beta}{2-\alpha -\beta}$.
The resulting equation has four roots when solved for $\beta$: the double root
 $\beta=1-\alpha$ and two roots $\beta=\pm \alpha$.
Solution
 $\beta=1-\alpha$ corresponds to the independent sequence with 
$r=0$. Since $\beta\geq 0$, therefore the only non-trivial solution is $\beta=\alpha$, which gives
$p=\frac12$ and  $X_k=\pm 1$.

Condition (\ref{(A)}) in this case is verified by
direct computation with conditional probabilities.
\end{proof}

\section{Auxiliary results and proofs}

Condition \rf{(A')} 
determines the form of the covariance matrix $r_k=E(X_0X_k)$. 
\begin{Lemma}\lbl{oldie1} Suppose that $(X_k)$ is an
$L_2$-stationary sequence such that
condition
\rf{(A')} holds true  and $2|r|<1+r_2$. 
Then $corr(X_0,X_k)=r^k$.
\end{Lemma}
\begin{proof} Indeed,  multiplying 
\rf{(A')} by $X_0$ we get
$r_k=a(r_{k-1}+r_{k+1})$.
In particular, if $r:=r_1 =0$ then $a=0$ and $r_k=0$ for all $k$. On the other hand, if
$r\ne 0$, then $1+r_2>0, a=r_1/(1+r_2)$ and 
the correlation coefficients $r_k$ satisfy the recurrence
$$ (1+r_2)r_k=r(r_{k-1}+r_{k+1}), k=1,2,\dots
$$
 From this we infer
that $r_k\to 0$ as $k\to \infty$. Indeed, since $|r_k|\leq 1$,
$r_\infty=\limsup_{k\to\infty}|r_k|$ is finite, and satisfies
$r_\infty(r_2+1)\leq 2r_\infty|r|$. 
Is is easy to see that since $r_k\to 0$, the recurrence has unique solution $r_k=r^k$.
\end{proof} 
We use the notation $E(\cdot|\dots,X_0)$ to denote the conditional
expectation with respect to the sigma field generated by $\{X_k:k\leq 0\}$.

The following Lemma comes from Bryc(1998); the proof is included for completeness.
\begin{Lemma}\lbl{oldie2}
If $(X_k)$ satisfies the assumptions of Lemma \ref{oldie1}, then
\be\lbl{LR0} 
E(X_1|\dots, X_0)=r X_0
\ee
\end{Lemma}
\begin{proof} By Lemma \ref{oldie1}, we have $r_k=r^k$, and $|r|<\frac{1+r^2}{2}\leq 1$.

We first show by induction that for all $n\in\ZZ, k\in\NN, 0\leq i\leq k$
\be \label{EE1}
E(X_{n+i}| \dots, X_{n-1}, X_n,X_{n+k},X_{n+k+1},\dots)=
a(i,k)X_n+b(i,k)X_{n+k} \hbox{}
\ee
where $a(i,k)=\frac{r_i-r_{k-i}r_k}{1-r_k^2}$,
 $b(i,k)=\frac{r_{k-i}-r_{i}r_k}{1-r_k^2}$

For $k=2$, (\ref{EE1}) follows from (\ref{(A')}) when $i=1$.  
Clearly,  (\ref{EE1}) trivially holds true
when $i=0$ or $i=k$ for all $k$. 

Suppose that
(\ref{EE1})  holds true for a given value of $k\geq 2$ and all 
$n\in\ZZ$. We will prove that it holds true for $k+1$. We  only need to show that the
left-hand side of \rf{EE1} is a linear function of the appropriate variables.
Indeed, in the non-degenerate case the
coefficients $a(i,k), b(i,k)$ in a linear regression are  uniquely determined from the
covariances; the covariance matrices are non-degenerate since $|r|<1$ and $r_k=r^k$.

Using routine properties of conditional expectations, 
the case of general index $0< i<k$ reduces to two values $i=1,k-1$. 
By symmetry, it suffices to give the
proof when  $i=1$.

Conditioning on additional variable
$X_{n+k}$ we get $$E(X_{n+1}|\dots, X_{n-1},
X_n,X_{n+k+1},X_{n+k+2},\dots)=$$ $$E(E^{\dots
X_{n-1}, X_n,X_{n+k},X_{n+k+1},\dots}(X_{n+1})| \dots
X_{n-1}, X_n,X_{n+k+1},\dots)=$$
$$a(1,k)X_n+b(1,k)E(X_{n+k}|\dots, X_{n-1}, X_n,X_{n+k+1},\dots)$$
Now adding $X_{n+1}$ 
to the condition we get.
$$E(X_{n+k}|\dots, X_{n-1},
X_{n},X_{n+k+1},X_{n+k+2},\dots)=$$ 
$$E(E^{\dots, X_{n},
X_{n+1},X_{n+k+1},X_{n+k+2},
\dots}(X_{n+k})|\dots, X_{n-1},
X_n,X_{n+k+1},X_{n+k+2},\dots)=$$       
$$a(k-1,k)E(X_{n+1}|\dots, X_{n-1},
X_n,X_{n+k+1},X_{n+k+2},\dots)+b(k,k)X_{n+k+1}$$
This gives the system of two linear equations for $E(X_{n+1}| \dots 
X_{n-1}, X_n,X_{n+k+1},X_{n+k+2},\dots)$, which has the unique solution which is a linear function of 
$X_n,X_{n+k+1}$
when $a(k-1,k)b(1,k)\ne 1$.
It remains to notice that if $k>1$ then $a(k-1,k)=b(1,k)=\frac{r^{k-1}-r^{k+1}}{1-r^{2k}}
<1$. Indeed, the latter is equivalent to $r^{k-1}(1-r+r^{k+1})<1$
and holds true because $-1<r<1$,  $r^{k+1}<1-r+r^{k+1}$, and $1-r+r^{k+1}\leq 1-r+r^{2}\leq 1$.

Therefore  the regression $E(X_{n+1}| \dots 
X_{n-1}, X_n,X_{n+k+1},X_{n+k+2},\dots)$ is linear, 
and (\ref{EE1})  holds for $k+1$.
This proves  (\ref{EE1}) by induction.

Passing to the limit as $k\to\infty$ in 
(\ref{EE1}) with $n=0, i=1$  we get \rf{LR0}.
\end{proof}

The following result comes from 
Bryc(1998).
Since certain minor details differ we include it here for completeness.
\begin{Lemma}\lbl{L 5.2} 
If $(X_k)$ satisfies the assumptions of Lemma \ref{oldie1} and \rf{(B')} holds true, then
\be \lbl{OldStuff}
(1-A(1+r^2))E(X_1^2|\dots, X_0)=(A(1-r^2)+Br^2)X_0^2 +C+D(1+r^2) X_0
\ee
\end{Lemma}
\begin{proof}
By Lemma \ref{oldie1} we have $L(x,y)=\frac{r}{1+r^2}(x+y)$. 
Since $E(X_1 X_2|{\dots, X_0})=E^{\dots, X_0}\left(X_1E^{\dots, X_1}(X_2)\right)$, from 
Lemma \ref{oldie2}
 we get 
\be\lbl{miss1}
E(X_1 X_2|{\dots, X_0})=rE(X_1^2|\dots, X_0)
\ee

We now give another expression for the left hand side of \rf{miss1}. Substituting
 $E(X_1 X_2|{\dots, X_0})=E(X_2E(X_1|\dots,X_0,X_2,\dots)|{\dots, X_0})=$ into (\ref{(A)})
we get $E(X_1 X_2|{\dots, X_0})= \frac{r}{1+r^2}E(X_2 (X_2+X_0)|{\dots, X_0})$. 
By Lemma \ref{oldie2}
this implies $E(X_1
X_2|{\dots, X_0})=\frac{r^3}{1+r^2}X_0^2+\frac{r}{1+r^2}E(X_2^2|{\dots, X_0})$.
Since $r\ne 0$, combining the latter with \rf{miss1} we have
\be \lbl{OldStuff2}
E(X_2^2|{\dots, X_0})=(1+r^2)E(X_1^2|{\dots, X_0})-r^2 X_0^2
\ee
We now substitute expression \rf{OldStuff2} in \rf{(B)} as follows. 
Taking the conditional expectation $E(\cdot|\dots,X_0)$ of both sides of
(\ref{(B)}), with $k=1$ and substituting \rf{Q},  we get
$$
E(X_1^2|{\dots, X_0})=A X_0^2+A E(X_2^2|{\dots, X_0})+B X_0^2 r^2+C +D(1+r^2) X_0
$$
Replacing $E(X_2^2|{\dots, X_0})$ by the right hand side of \rf{OldStuff2}
 we get \rf{OldStuff}.
\end{proof}

The following result serves as a lemma but is of independent interest.

\begin{Proposition}\lbl{T1} Suppose $(X_k)$ is a square-integrable standardized
stationary homogeneous Markov chain such that 
 the correlation coefficients  satisfy $r\ne 0,2|r|<1+r_2$.

If $(X_k)$ satisfies condition
(\ref{(A)}) and the conditional variance $Var(X_k|X_{k-1})$ is
 a quadratic function of $X_{k-1}$ then
one of the following condition holds true:

\be \lbl{const}
Var(X_k|X_{k-1})=const
\ee
or
\be \lbl{Qdr}
Var(X_k|X_{k-1})=(1-r^2)X_{k-1}^2
\ee
\end{Proposition}

\begin{Remark}\lbl{Siva} Condition \rf{Qdr} implies that
$|X_k|=|X_{k-1}|$ for all $k$, even in the non-Markov case. 
\end{Remark}

\begin{Remark}
If linear regression condition (\ref{(A)}) is weakened to a symmetric pair
of conditions $E(X_k|X_{k-1})=r X_{k-1}$ and  $
E(X_{k-1}|X_{k})=r X_{k}$ then the conditional
variance can be given by other quadratic expressions, see Example \ref{Chebyshev}.
\end{Remark}

\begin{proofof}{Proposition \ref{T1}} 
If $Var(X_k|X_{k-1})$ is quadratic then there are constants $a,b,c$ such that

\be \lbl{*}
E(X_k^2|X_{k-1})=a X_{k-1}^2 + b X_{k-1} +c
\ee

Since $(X_k)$ is a homogeneous Markov chain and \rf{LR0} holds true
\be \lbl{**}
E(X_{k+1}^2|X_{k-1})=E(a X_k^2+bX_k+c|X_{k-1})=a^2 X_{k-1}^2+(a+r)bX_{k-1}+(a+1)c
\ee

On the other hand, condition (\ref{(A)}) implies, see (\ref{OldStuff2})
\be
(1+r^2)E(X_k^2|X_{k-1})=r^2X_{k-1}^2+E(X_{k+1}^2|X_{k-1})
\ee

Combining this  with (\ref{*}) and (\ref{**}) we get
\be \lbl{EqtnSystem}
(1+r^2)aX_{k-1}^2+(1+r^2)bX_{k-1}+(1+r^2)c=(a^2 +r^2)X_{k-1}^2+(a+r)bX_{k-1}+(a+1)c
\ee

Since $E(X_{k-1})=0$ and $E(X_{k-1}^2)=1$ therefore $X_{k-1}$ must have at least two values. 
We consider separately two cases.
\begin{itemize}
\item[(a)]  If $X_k$ has only two values then by Proposition \ref{L1}
 $X_k=\pm 1$   and $Var(X_k|X_{k-1})=1-r^2$ is a non-random constant, ending the proof. 

\item[(b)] If $X_{k-1}$ has at least three values, then $X_{k-1}^2,X_{k-1}, 1$ are linearly
independent. Therefore (\ref{EqtnSystem}) implies 
\be \lbl{System2}
(1+r^2)a=a^2+r^2,
(1+r^2)b=(a+r)b,
(1+r^2)c=(a+1)c
\ee 
Since \rf{*} implies that $a+c=1$, the only solutions of (\ref{System2}) are  $c\ne 0, a=r^2$ or  $c=0,a=1$.
Since $0<|r|<1$, both solutions imply $b=0$. 

Clearly,
$a=r^2$ implies \rf{const}.
On the other hand if $c=0$ and $a=1$,
then  $E(X_k^2|X_{k-1})=X_{k-1}^2$.  Thus \rf{Qdr} hold true.
\end{itemize}
\end{proofof}


\begin{proofof}{Theorem \ref{Corollary}}
We first consider the two-valued case. 
If $X_{k-1}^2$ is a non-random constant, then $X_{k-1}^2=1$ and 
thus $Q$ is non-unique; one can take
$Q(x,y)=(x^2+y^2)/2$ to satisfy (\ref{degenerate1}), or one can take
$Q(x,y)=\frac{r^2}{1+r^4}(x^2+y^2)+\frac{(1-r^2)^2}{1+r^4}$ to satisfy
(\ref{MainAssumption}). 

Suppose now that $X_k$ has more than two values.
We first verify that that the collusion \rf{degenerate1}  holds true when
$A=1/(1+r^2)$. 
In this case the left hand side of (\ref{OldStuff}) is zero.
Since $X_k$ has more than two values, this implies that
$D=0$ and $C=0$. Therefore (\ref{C}) implies (\ref{degenerate1}).

Now consider the case when $A\ne 1/(1+r^2)$.
From (\ref{OldStuff}) we have
\be \lbl{Var}
E(X_k^2|X_{k-1})=\frac{A(1-r^2)+Br^2}{1-A(1+r^2)} X_{k-1}^2+\alpha X_{k-1}+\beta 
\ee
where $\alpha=\frac{D(1+r)}{1-A(1+r^2)}$.
This shows that  $Var(X_k|X_{k-1})$
is quadratic. 
%
By Proposition \ref{T1} we have $\alpha=0$; since $|r|<1$ this implies that $D=0$.
 We also know that either 
\rf{const} holds true,
which is equivalent to $E(X_k^2|X_{k-1})=r^2 X_{k-1}^2+1-r^2$,
or \rf{Qdr} holds true, which is equivalent to $E(X_k^2|X_{k-1})=X_{k-1}^2$. We now compare these
two expressions with (\ref{Var}): since $\alpha=0$ and $X_{k-1}^2$ is non-constant,  the coefficients at $X_{k-1}^2$
must match. That is, either  $\frac{A(1-r^2)+Br^2}{1-A(1+r^2)}=r^2$
or $\frac{A(1-r^2)+Br^2}{1-A(1+r^2)}=1$.
By a simple algebra the former implies (\ref{MainAssumption})
and the latter implies (\ref{degenerate1}).
\end{proofof}


\begin{Lemma} \lbl{LemmaBounds}
Suppose that $E(X)=E(Y)=0, E(X^2)=E (Y^2)=1, E(X^4)= E(Y^4)<\infty$ 
and the following conditions hold true
\begin{itemize}
\item $E(Y|X)=r X$  
\item $E(X^3|Y)=\alpha Y^3 +\beta Y$
\item $\alpha \ne r$
\end{itemize}
Then $\frac{\beta}{r-\alpha}\geq 1$.

\end{Lemma}
\begin{proof}
Conditioning in two different directions in $E X^3 Y$ we get
$r EX^4= \alpha E(Y^4)+\beta E(Y^2)$.
Therefore  $E (X^4)= \frac{\beta}{r-\alpha}$.
Since $ E( X^4)\geq (E(X^2))^2 = 1 $ we have  $\frac{\beta}{r-\alpha}\geq 1$, which
ends the proof.
\end{proof}

The following lemma is based on 
estimates from Bryc(1995), Theorem 6.2.2.
The proof is omitted.
\begin{Lemma}\lbl{B-P}
Suppose $X,Y$ are square-integrable random variables 
with the same distribution. Let $r=corr(X,Y)$ denote the correlation coefficient and assume that
$r\ne 0,\pm 1, E(X|Y)=rY, E(Y|X)=rX, Var(X|Y)=1-r^2, Var(Y|X)=1-r^2$. 
Then $E(X^4)\leq  32\frac{r^2+2|r|+2}{(1-|r|)r^4} $.
\end{Lemma}

%
%

\begin{proofof}{Theorem \ref{T4}}
Since the conclusion is trivially true when $A= 1/(1+r^2)$,
throughout the proof we  assume that
$A\ne 1/(1+r^2)$. 
In this case \rf{OldStuff} implies  $Var(X_k|X_{k-1})=1-r^2$. Since the assumptions
are symmetric, and $0<|r|<1$ by Lemma \ref{B-P} and stationarity
 we have $E(X_1^4)=E(X_2^4)<\infty$.

Notice that \rf{OldStuff} implies 
$E(X_2^2|X_0)=E^{X_0}(E(X_2^2|\dots, X_{1})=r^2E(X_1^2|{X_0})+1-r^2$.
Thus
\be\lbl{miss2}
E(X_2^2|X_0)=
r^4X_0^2+1-r^4
\ee

We now compute conditional moments using the approach of Plucinska(1983). 
Using constant conditional variance and
 (\ref{(A')}), we  write $E(X_1 X_2^2|X_0)$ in two different
ways  as $$E(E(X_1 X_2^2|\dots, X_0, X_1)|X_0)=E(r^2X_1^3+(1-r^2)X_1|X_0)$$ and
 as 
$$E(E(X_1 X_2^2|X_2, X_0)|X_0)=\frac{r}{1+r^2}E(X_2^2(X_2+X_0)|X_0)$$
Combining these two representations and using \rf{miss2}, and $r\ne0$
 we get after simple algebra
\be\lbl{eq1}
r E(X_1^3|X_0)=\frac{1}{1+r^2}E(X_2^3|X_0)+\frac{r^4}{1+r^2} X_0^3
\ee
Similarly, we rewrite $E(X_1^2 X_2|X_0)$ in two different
ways as $$E(E(X_1^2 X_2|\dots, X_0, X_1)|X_0)=rE(X_1^3|X_0)$$ and, using (\ref{(B')}), as 
$$E(E(X_1^2 X_2|X_2, X_0)|X_0)=E((A(X_2^2+X_0^2)+BX_0X_2+C) X_2|X_0)$$
Using \rf{miss2}, after some algebra we get
\be\lbl{eq2}
r E(X_1^3|X_0)=r^2(A+B r^2) X_0^3 + A E(X_2^3|X_0)+(B(1-r^4)+Cr^2)X_0 
\ee

Solving  the system of equations \rf{eq1}, \rf{eq2} for $E(X_1^3|X_0)$ we get
\be
E(X_1^3|X_0)=r\frac{A(1-r^2)+Br^2}{1-A(1+r^2)}X_0^3+\frac{B(1-r^4)+Cr^2}{r(1-A(1+r^2))}X_0
\ee
Substituting
\rf{C},  \rf{MainAssumption}, and  denoting $\tilde{A}=A(1+r^2)$ we have

\be
E(X_1^3|X_0)=r^3 X_0^3 - \frac{1-r^2}{r^3}
\frac{\tilde{A}(1+2r^4)-r^2(1+2r^2)}{1-\tilde{A}}X_0
\ee

Therefore by Lemma \ref{LemmaBounds} and a simple calculation we have
\be
\frac{\tilde{A}(1+r^4)-r^2(1+r^2)}{r^4(1-\tilde{A})}\leq 0
\ee
Since $\frac{r^2}{1+r^4}<\frac{1}{1+r^2}$ this implies that either $A>1/(1+r^2)$ or
$
A\leq \frac{r^2}{1+r^4}
$.
\end{proofof}

%

\begin{proofof}{Theorem \ref{T+}}
 For $A\ne 1/(1+r^2)$ let 
\be \lbl{q}
q=\frac{r^2-A(1+r^2)}{r^4(1-A(1+r^2))}
\ee
The range of values of $A$ implies that $-1\leq q\leq 1$.
We give the proof for the case $-1< q \leq 1$. 
The only change needed for the case  $q=-1$, is
 to use  the symmetric two-valued Markov chain
defined in Section \ref{TwoValued} instead of the Markov chain $M_k$ defined below.

Define orthogonal polynomials $Q_n(x)$ by the recurrence
\be \lbl{q-recurrence}
Q_{n+1}(x)=x Q_n(x) - (1+q+\dots+q^{n-1})Q_{n-1}(X)
\ee
with $Q_0(x)=1$, $Q_1(x)=x$.
Let $\mu(dx)$ denote the probability measure which orthogonalizes $Q_n$
(see eg Chihara(1978), Theorem 6.4),
 and for fixed 
$-1<r<1$ define
\be \lbl{M-transition}
P(x,dy)=\sum_{n=0}^\infty r^n \tilde{Q}_n(x)\tilde{Q}_n(y)\mu(dy)
\ee 
where $\tilde{Q}_n(x)=Q_n(x)/\|Q_n\|_{L_2(\mu)}$ are normalized orthogonal polynomials
$Q_n$.
By Bryc(1998), Lemma 8.1,  for $-1<q\leq 1$ formula 
(\ref{M-transition}) defines a Markov transition
function with invariant measure
$\mu$. 
For $-1<q\leq 1$, let $M_k$ be a stationary Markov chain with the initial 
distribution $\mu$ and transition probability $P(x,dy)$.

It is known that $\mu$ is either gaussian or of bounded support, see
Koekoek-Swarttouw(1994), and hence the joint distribution of $M_1,\dots,M_d$ is
uniquely determined by mixed moments
$E(M_1^{k_1}\dots,M_d^{k_d})$.
 We will show by induction
with respect to $d$ that
\be \lbl{induct}
E(X_1^{k_1}\dots, X_d^{k_d})=E(M_1^{k_1}\dots M_d^{k_d})
\ee for all $d\geq 1$ and all non-negative integers $k_1,...,k_d$.

By Bryc(1998) marginal distributions are equal, $X_1\conj M_1$;
this shows that equality (\ref{induct}) holds true for all integer $k_1\geq 0$ when
$d=1$. Suppose \rf{induct} holds for all $k_1,...,k_d\geq 0$. Fix integer $k=k_{d+1}\geq 0$.
Expand polynomial $x^k$ into orthogonal expansion,
 $x^k=\sum_{j=0}^k a_jQ_j(x)$. Then 

$$E(X_1^{k_1}\dots, X_d^{k_d}X_{d+1}^k)= 
\sum_{j=0}^k a_j E(X_1^{k_1}\dots, X_d^{k_d}E(Q_j(X_{d+1})|X_1,\dots,X_{d})$$

Repeating the reasoning that lead to Bryc(1998), Lemma 6.3, we have
$E(Q_j(X_{d+1})|X_1,\dots,X_{d})=r^j Q_{j}(X_d)$.
Therefore $E(X_1^{k_1}\dots  X_d^{k_d}X_{d+1}^k)=\sum r^j a_j
E(X_1^{k_1}\dots X_d^{k_d}Q_j(X_{d}))$ is expressed as a linear
combination of moments that involve only $E(X_1^{j_1}\dots  X_d^{j_d})$.
Since
the
same reasoning applies to $M_k$, we have 
$E(M_1^{k_1}\dots M_d^{k_d}M_{d+1}^k)=\sum r^j a_j
E(M_1^{k_1}\dots M_d^{k_d}Q_j(M_{d}))$, and
(\ref{induct}) follows.   
\end{proofof}
\section{Example}

This section contains an example
of a stationary reversible Markov chain with linear regressions and quadratic
conditional moments, which does not satisfy condition (\ref{(A)}). The Markov chain has polynomial
regressions of all orders, and does not satisfy the conclusion of
Proposition \ref{T1}.

\begin{Example}\lbl{Chebyshev}
Suppose $T_n(x)$ are Chebyshev polynomials of the first kind,
$T_0=1, T_1(x)=x, T_2(x)=2x^2-1, x T_n(x)=\frac12 T_{n+1}(x)+\frac{1}{2} T_{n-1}(x)$.
Let $\mu(dx)=\frac{1}{\pi}\frac{1}{\sqrt{1-x^2}} dx$. 
Then $T_n$ are orthogonal in $L_2(d\mu)$ and
$\|T_0\|_{L_2(d\mu)}^2=1$ and for $k>0$ $\|T_k\|_{L_2(d\mu)^2}^2=\frac{1}{2}$.
We  define transition density by $p(x,y)=\sum_{n=0}^\infty r^n T_n(x)T_n(y)$.

Since $T_n(x)=\cos(n \arccos(x))$, the series can be summed.
Writing $T_n(x)=\cos(n \theta_x)$ we have $T_n(x) T_n(y)=
\frac12 \cos(n (\theta_x+\theta_y))+\frac12 \cos(n 
(\theta_x-\theta_y))$ 
Therefore $$p(x,y)=\frac12
 \frac{1-r\cos(\theta_x+\theta_y)}{1+r^2-2r \cos(\theta_x+\theta_y)}
+\frac12
 \frac{1-r\cos(\theta_x-\theta_y)}{1+r^2-2r\cos(\theta_x-\theta_y)}$$
This shows that $p(x,y)\geq \frac{1-|r|}{(1+|r|)^2}>0$. The expression simplifies to

$$p(x,y)=
\frac{1-r^2+ r(2 r (y^2+y^2)-(3 +r^2) y x)}
{(1-r^2)^2+4r^2 (x^2+y^2- (r+1/r)x y)}
$$

Thus we can define the Markov chain $X_k$ with one-step transition probabilities
$P_x(dy)=p(x,y)\mu(dy)$ and initial distribution $\mu$. 
Since $\int p(x,y)\mu(dx)=1$, the chain is stationary.

Notice that by the definition of $p(x,y)$ we have 
$E(T_n(X_1)|X_0)=r^n \|T_n\|_2^2  T_n(X_0)$. 
Therefore for $n\ge 1$ we have $E(T_n(X_1)|X_0)=\frac12 r^nT_n(X_0)$ 

In particular
$E(X_1|X_0)=r/2 X_0$, and  
$E(2X_1^2-1|X_0)=\frac12r^2 (2X_0^2-1)$. The latter implies 
$E(X_1^2|X_0)= \frac12r^2 X_0^2+\frac12-\frac14r^2$ and hence the conditional variance
$Var(X_1|X_0)=\frac14r^2 X_0^2+\frac12-\frac14r^2$ is non-constant.
This should be contrasted with the conclusion of Proposition \ref{T1} 
and assumptions in
{Bryc(1998), Wesolowski(1993)}.
\end{Example}

{\bf Acknowledgements} I would like to thank Dr. Sivaganesan and W. Matysiak 
for helpful discussions.

\end{document}